\documentclass{au}

\usepackage{latexcad}
\usepackage{natbib}

\def\ldiv{\setminus}
\def\rdiv{/}

\theoremstyle{definition}
\newtheorem{definition}{Definition}[section]
\newtheorem{example}[definition]{Example}
\newtheorem{remark}[definition]{Remark}

\theoremstyle{theorem}
\newtheorem{lemma}[definition]{Lemma}

\newtheorem{proposition}[definition]{Proposition}

\title[Loops of Bol-Moufang Type]{The Varieties of Loops of Bol-Moufang Type}

\author{J.~D.~Phillips}

\address{Department of Mathematics \& Computer Science, Wabash College,
Crawfordsville, Indiana 47933, U.S.A.}

\email{phillipj@wabash.edu}

\author{Petr Vojt\v{e}chovsk\'y}

\address{Department of Mathematics, University of Denver, 2360 S Gaylord St,
Denver, Colorado 80208, U.S.A.}

\email{petr@math.du.edu}

\begin{document}

\begin{abstract}
A loop identity is of Bol-Moufang type if two of its three variables occur once
on each side, the third variable occurs twice on each side, and the order in
which the variables appear on both sides is the same, viz. $((xy)x)z=x(y(xz))$.
Loop varieties defined by one identity of Bol-Moufang type include groups, Bol
loops, Moufang loops and C-loops. We show that there are exactly $14$ such
varieties, and determine all inclusions between them, providing all necessary
counterexamples, too. This extends and completes the programme of Fenyves
\cite{Fe2}.
\end{abstract}



\maketitle

\section{Introduction}

\noindent An identity $\varphi=\psi$ is said to be of \emph{Bol-Moufang type}
if: (i) the only operation appearing in $\varphi=\psi$ is a binary operation,
(ii) the number of distinct variables appearing in $\varphi$ (and thus in
$\psi$) is $3$, (iii) the number of variables appearing in $\varphi$ (and
thus in $\psi$) is $4$, (iv) the order in which the variables appear in
$\varphi$ coincides with the order in which they appear in $\psi$.

The \emph{variety of loops} consists of universal algebras
$(L,\,\cdot,\,\ldiv,\,\rdiv,\,e)$ whose binary operations $\cdot$, $/$,
$\setminus$ satisfy
\begin{displaymath}
    a\cdot(a\ldiv b) = b,\quad (b\rdiv a)\cdot a = b,\quad a\ldiv(a\cdot b) = b,
    \quad (b\cdot a)\rdiv a = b,
\end{displaymath}
and whose nullary operation $e$ satisfies
\begin{displaymath}
    e\cdot a = a\cdot e = a.
\end{displaymath}
Recall that any loop satisfies the identities $(x\rdiv y)\ldiv x =y$,
$x\rdiv(y\ldiv x)=y$.

For the rest of the paper, all identities of Bol-Moufang type $\mathcal B$ are
loop identities with $\cdot$ as the binary operation.

We say that all identities in a subset of $\mathcal B$ are \emph{equivalent}
if each of them defines the same variety of loops. In this sense, most of the
varieties of Bol-Moufang type can be defined in several equivalent ways. It
is then a nontrivial task of practical importance to describe all maximal
subsets of equivalent identities in $\mathcal B$. In fact, this work was
partially motivated by the authors' frustration with the inconsistencies in
the literature concerning definitions of loop varieties.

This paper presents the classification of all varieties of \emph{loops of
Bol-Moufang type}, i.e., varieties of loops defined by one identity from
$\mathcal B$. We determine the variety defined by each identity $\mathcal B$
and conclude that there are 14 such varieties, including groups, Bol loops,
Moufang loops and C-loops. We then describe the inclusions among all these
varieties and provide all necessary distinguishing examples.

Many of the results below were known already to Fenyves \cite{Fe1},
\cite{Fe2}. See the Acknowledgement and historical remarks for more
information.

\section{Systematic notation}

\noindent Let $x$, $y$, $z$ be all the variables appearing in the identities of
$\mathcal B$. Without loss of generality, we can assume that they appear in the
terms in alphabetical order. Then there are exactly $6$ ways in which the $3$
variables can form a word of length $4$, and there are exactly $5$ ways in
which a word of length $4$ can be bracketed, namely:
\begin{displaymath}
\begin{array}{cc}
    \begin{array}{c|c}
        A&xxyz\\
        B&xyxz\\
        C&xyyz\\
        D&xyzx\\
        E&xyzy\\
        F&xyzz
    \end{array}
    \quad\quad&
    \begin{array}{c|c}
        1&o(o(oo))\\
        2&o((oo)o)\\
        3&(oo)(oo)\\
        4&(o(oo))o\\
        5&((oo)o)o
    \end{array}
\end{array}
\end{displaymath}
Let $Xij$ with $X\in\{A,\dots,F\}$, $1\le i<j\le 5$ be the identity from
$\mathcal B$ whose variables are ordered according to $X$, whose left-hand side
is bracketed according to $i$, and whose right-hand side is bracketed according
to $j$. For instance, $C25$ is the identity $x((yy)z)=((xy)y)z$.

It is now clear that any identity in $\mathcal B$ can be transformed into some
identity $Xij$ by renaming the variables and interchanging the left-hand side
with the right-hand side. There are therefore $6\cdot(4+3+2+1)=60$
``different'' identities in $\mathcal B$.

The \emph{dual} of an identity $I$ is the identity obtained from $I$ by reading
it backwards, i.e., from right to left. For instance, the dual of
$(xy)(xz)=((xy)x)z$ is the identity $z(x(yx)) = (zx)(yx)$. With the above
conventions in mind, we can rewrite the latter identity as $x(y(zy)) =
(xy)(zy)$. One can therefore identify the dual of any identity $Xij$ with some
identity $X'j'i'$. The name $X'j'i'$ of the dual of $Xij$ is easily calculated
with the help of the following rules:
\begin{displaymath}
    A'=F,\quad B'=E,\quad C'=C,\quad D'=D,\quad 1'=5,\quad 2'=4,\quad 3'=3.
\end{displaymath}

Finally, we will use the following notational conventions: we usually omit
$\cdot$ while multiplying two elements (eg $x\cdot y = xy$), we declare $\ldiv$
and $\rdiv$ to be less binding than the omitted multiplication (eg $x\rdiv yz =
x\rdiv (yz)$), and if $\cdot$ is used, we consider it to be less binding than
any other operation (eg $x\cdot yz\ldiv y = x((yz)\ldiv y)$).

\section{Canonical definitions of some varieties of loops}

\noindent Table \ref{Tb:Definitions} defines $15$ varieties of loops. With the
exception of the $3$-power associative loops, all these varieties can be
defined by some identity $Xij$. Namely, GR is equivalent to $A12$ (cancel $x$
on the left), LA to $A45$ (substitute $e$ for $z$), RA to $F12$ (duality), and
FL to $B45$ (substitute $e$ for $z$).

\begin{table}\caption{Definitions of varieties of loops of Bol-Moufang type.}\label{Tb:Definitions}
\begin{displaymath}
    \begin{array}{ccccc}
        \text{variety}&\text{abbrev.}&\text{defining
        identity}&\text{its name}&\text{ref.}\\
        \hline
        \text{groups}&
            \text{GR}&x(yz)=(xy)z&&\text{folklore}\\
        \text{extra loops}&
            \text{EL}&x(y(zx))=((xy)z)x&D15&\text{\cite{Fe1}}\\
        \text{Moufang loops}&
            \text{ML}&(xy)(zx)=(x(yz))x&D34&\text{\cite{Mo},
            \cite[p.\ 58]{Br}, \cite{Pf}}\\
        \text{left Bol loops}&
            \text{LB}&x(y(xz))=(x(yx))z&B14&\text{\cite{Ro}}\\
        \text{right Bol loops}&
            \text{RB}&x((yz)y)=((xy)z)y&E25&\text{\cite[p.\ 116]{Br},
            \cite{Ro}}\\
        \text{C-loops}&
            \text{CL}&x(y(yz))=((xy)y)z&C15&\text{\cite{Fe2}}\\
        \text{LC-loops}&
            \text{LC}&(xx)(yz)=(x(xy))z&A34&\text{\cite{Fe2}}\\
        \text{RC-loops}&
            \text{RC}&x((yz)z)=(xy)(zz)&F23&\text{\cite{Fe2}}\\
        \text{left alternative loops}&
            \text{LA}&x(xy)=(xx)y&&\text{folklore}\\
        \text{right alternative loops}&
            \text{RA}&x(yy)=(xy)y&&\text{folklore}\\
        \text{flexible loops}&
            \text{FL}&x(yx)=(xy)x&&\text{\cite[p.\ 89]{Pf}}\\
        \text{left nuclear square loops}&
            \text{LN}&(xx)(yz)=((xx)y)z&A35&\\
        \text{middle nuclear square loops}&
            \text{MN}&x((yy)z)=(x(yy))z&C24&\\
        \text{right nuclear square loops}&
            \text{RN}&x(y(zz))=(xy)(zz)&F13&\\
        \text{$3$-power associative loops}&
            \text{$3$PA}&x(xx)=(xx)x&&
    \end{array}
\end{displaymath}
\end{table}

We have carefully chosen the defining identities in such a way that they are
either self-dual (GR, EL, CL, FL, MN, $3$PA) or coupled into dual pairs
(LB$'$=RB, LC$'$=RC, LA$'$=RA, LN$'$=RN). The only exception to this rule is
the Moufang identity $D34$. We will often appeal to this duality in our proofs.

Only four of the above varieties were not previously named in the literature,
namely the \emph{left}, \emph{middle} and \emph{right nuclear square loops},
and the \emph{$3$-power associative loops}. Since we have no desire to swamp
the field with new definitions, we opted for these longer, descriptive names.
The reader should note that $3$-power associative loops are not necessarily
\emph{power associative}, i.e., the subloop generated by $x$ does not have to
be a group in a $3$-power associative loop (viz. Example \ref{Ex:3PAnPA}).
We have included $3$PA as a technical variety that will allow us to make
several arguments faster.

\begin{example}\label{Ex:3PAnPA}
This is a loop that is $3$-power associative but not power
associative (since $(1\cdot 1)(1\cdot 1)\ne 1(1(1\cdot 1))$):
\begin{displaymath}\begin{array}{cccccc}
    0&1&2&3&4&5\\
    1&2&0&4&5&3\\
    2&0&3&5&1&4\\
    3&4&5&0&2&1\\
    4&5&1&2&3&0\\
    5&3&4&1&0&2
\end{array}\end{displaymath}
\end{example}

A loop $L$ is said to have the \emph{left inverse property} if $e\rdiv x\cdot
xy = y$ for every $y\in L$. Dually, $L$ has the \emph{right inverse property}
if $yx\cdot x\ldiv e=y$ for every $y\in L$. If $L$ has both the left and
right inverse property, it is called an \emph{inverse property loop}.

When $L$ has the left inverse property, it also has two-sided inverses since
$e\rdiv x = e\rdiv x\cdot x(x\ldiv e) = x\ldiv e$. The same conclusion holds
when $L$ has the right inverse property. The two-sided inverse of $x$ will be
denoted by $x^{-1}$. However, when a loop has two-sided inverses, it does not
have to have any inverse properties, as the following example shows:

\begin{example}\label{Ex:nLIP}
This is a loop that has two-sided inverses but is neither a left inverse
property loop (since $1^{-1}(1\cdot 2)\ne 2$), nor a right inverse property
loop (since $(2\cdot 1)1^{-1}\ne 2$):
\begin{displaymath}\begin{array}{ccccc}
    0&1&2&3&4\\
    1&0&3&4&2\\
    2&4&0&1&3\\
    3&2&4&0&1\\
    4&3&1&2&0
\end{array}\end{displaymath}
\end{example}

\section{Equivalences}\label{Sc:E}

\noindent We now begin the exhaustive search for equivalent Bol-Moufang
identities.

\begin{proposition}\label{Pr:First}
The following Bol-Moufang identities are equivalent to
the defining group identity $x(yz)=(xy)z$: $A12$, $A23$, $A24$, $A25$, $B12$,
$B13$, $B24$, $B25$, $B34$, $B35$, $C13$, $C23$, $C34$, $C35$, $D12$, $D13$,
$D14$, $D25$, $D35$, $D45$, $E13$, $E14$, $E23$, $E24$, $E35$, $E45$, $F14$,
$F24$, $F34$, $F45$.
\end{proposition}
\begin{proof} We have already noted that $A12$ defines groups. We now briefly
describe how each of the remaining identities listed in this Proposition can be
seen to be equivalent to groups.

For $A23$: let $z=e$, deduce LA, then use LA to rewrite $A23$ into $A12$. For
$A24$: note that given $x$, any $u$ can be written as $xy$ for some $y$; then
use $u=xy$ in $A24$. For $A25$: let $z=e$, deduce LA, then use LA to rewrite
$A25$ into $A24$. For $B12$: cancel $x$ on the left. For $B13$: let $u=xz$.
For $B24$: let $u=yx$. For $B25$: let $z=e$, deduce FL, then use FL to
rewrite $B25$ as $B24$. For $B34$: let $z=e$, deduce FL, then use FL to write
$B34$ as $B35$, let $u=xy$ in $B35$. For $B35$: see $B34$. For $C13$: let
$u=yz$. For $C23$: let $x=e$, deduce LA, then use LA to write $C23$ as $C13$.
For $D12$: cancel $x$ on the left. For $D13$: let $u=zx$. For $D14$: let
$z=e$, deduce FL, then use FL to write $D14$ as $D12$.

The remaining $15$ identities are duals of the already investigated $15$
identities. Since the defining group identity is self-dual, we are done.
\end{proof}

\begin{proposition}\label{Pr:Extra}
The following Bol-Moufang identities are equivalent to
the defining extra identity $D15$: $B23$, $D15$, $E34$.
\end{proposition}
\begin{proof} See \cite[Thm.\ 1]{Fe1}.
\end{proof}

\begin{proposition} The following Bol-Moufang identities are equivalent to
the defining Moufang loop identity $D34$: $B15$, $D23$, $D34$, $E15$.
\end{proposition}
\begin{proof}
In \cite[p.88--89]{Pf}, Pflugfelder defines identities ($MI$), ($M_5$), ($M_6$)
and $(M_7)$, and shows in \cite[Thm.\ IV.1.4]{Pf} that these identities are
equivalent. They correspond to our identities $D34$, $D23$, $B15$ and $E15$,
respectively.
\end{proof}

\begin{remark} We would like to point out that \cite[Thm.\ IV.1.4]{Pf} also
says that ($MI$) is equivalent to ($M_4$)=B14. This is true only if flexibility
holds in the loop in question, which is what Pflugfelder tacitly assumes.
\end{remark}

\begin{lemma}\label{Lm:LC1}
Let $L$ be an LC-loop. Then $L$ is left alternative, has the left inverse
property, is a middle nuclear square loop, and satisfies $C14$.
\end{lemma}
\begin{proof}
The left alternative law follows from $A34$ with $z=e$. Hence $A34$ implies
$A14$. By $A14$, we have $x(x\cdot (x\ldiv e)z)=x(x\cdot x\ldiv e)\cdot
z=xz$, and thus $x\cdot (x\ldiv e)z=z$. With $x=e\rdiv y$, we obtain $z =
e\rdiv y\cdot ((e\rdiv y)\ldiv e)z = e\rdiv y\cdot yz$, and $L$ has the left
inverse property.

By $A14$ and the left inverse property, $x(x\cdot (x^{-1}y)z) = x(x\cdot
x^{-1}y)\cdot z = xy\cdot z$. With $(x^{-1}y)^{-1}z$ instead of $z$, we get
$xy\cdot (x^{-1}y)^{-1}z = x(xz)=(xx)z$. Therefore $(x^{-1}y)^{-1}z =
(xy)^{-1}\cdot (xx)z$, which reduces to $(x^{-1}y)^{-1} = (xy)^{-1}(xx)$ with
$z=e$. But then $(xy)^{-1}\cdot (xx)z = (x^{-1}y)^{-1}z = (xy)^{-1}(xx)\cdot
z$, and thus $L$ is a middle nuclear square loop. The identity $C14$ follows
by LA.
\end{proof}

\begin{lemma}\label{Lm:LC2}
Assume that $L$ is a loop satisfying $C14$. Then $L$ is an LC-loop.
\end{lemma}
\begin{proof}
The left alternative law follows from $C14$ with $x=e$. By $C14$, we have
$e=x\cdot x\ldiv e = (x\rdiv xx\cdot xx)(x\ldiv e) = (x\rdiv xx)(x\cdot
x(x\ldiv e)) = x\rdiv xx\cdot x$, and hence $e\rdiv x = x\rdiv xx$, or
$e\rdiv x\cdot xx = x$. Then $C14$ yields $e\rdiv x\cdot xy = e\rdiv x\cdot
x(x\cdot x\ldiv y) = (e\rdiv x\cdot xx)(x\ldiv y) = x(x\ldiv y) = y$, and $L$
has the left inverse property.

By $C14$ and the left inverse property, $x(yy)\cdot y^{-1}(y^{-1}x^{-1}) =
x\cdot y(y\cdot y^{-1}(y^{-1}x^{-1})) = e$, and thus $(x\cdot
yy)^{-1}=y^{-1}\cdot y^{-1}x^{-1}$. Applying the left inverse property to
$(x\cdot yy)z = x(y\cdot yz)$ yields $y^{-1}(y^{-1}(x^{-1}\cdot (x\cdot
yy)z)) = z$. With $u=(x\cdot yy)z$, the last identity becomes
$y^{-1}(y^{-1}(x^{-1}u)) = z=(x\cdot yy)^{-1}u =(y^{-1}\cdot y^{-1}x^{-1})u$,
and $A34$ follows by LA.
\end{proof}

\begin{proposition} The following Bol-Moufang identities are equivalent to
the defining LC-loop identity $A34$: $A14$, $A15$, $A34$, $C14$.
\end{proposition}
\begin{proof} Note that any of the three identities $A14$, $A15$, $A34$ yield
LA (let $y=e$ in $A14$, $z=e$ in $A15$, use Lemma \ref{Lm:LC1}(i) for $A34$).
With LA, the three identities are immediately seen to be equivalent. Lemmas
\ref{Lm:LC1} and \ref{Lm:LC2} show that C14 is equivalent to A34.
\end{proof}

By the duality, we obtain:

\begin{proposition} The following Bol-Moufang identities are equivalent to
the defining RC-loop identity $F23$: $C25$, $F15$, $F23$, $F25$.
\end{proposition}

\begin{proposition} The following Bol-Moufang identities are equivalent to
the defining left alternative identity $x(xy)=(xx)y$: $A13$, $A45$, $C12$.
\end{proposition}
\begin{proof} Both $A13$ and $A45$ yield LA with $z=e$, while $C12$ yields LA
with $x=e$. On the other hand, LA obviously implies each of the three
identities. They are therefore equivalent.
\end{proof}

By the duality, we obtain:

\begin{proposition} The following Bol-Moufang identities are equivalent to
the defining right-alternative identity $x(yy)=(xy)y$: $C45$, $F12$, $F35$.
\end{proposition}

\begin{proposition}\label{Pr:Last}
The following Bol-Moufang identities are equivalent to
the defining flexible identity $x(yx)=(xy)x$: $B45$, $D24$, $E12$.
\end{proposition}
\begin{proof} As the defining identity FL is self-dual and $(B45)'=E12$, it
suffices to show that the identities $B45$ and $D24$ are equivalent to FL. It
is obviously true for $B45$ (cancel $z$ on the right). With $z=e$, $D24$
reduces to FL. Using FL with $yz$ instead of $y$ yields $D24$.
\end{proof}

The only identities not covered by Propositions \ref{Pr:First}--\ref{Pr:Last}
are  $A35$, $B14$, $C15$, $C24$, $E25$ and $F13$. These are the defining
identities of LN, LB, CL, MN, RB and RN, respectively.

\section{Implications}\label{Sc:I}

%FIGURE
\setlength{\unitlength}{0.8mm}
\begin{figure}[t]\centering\begin{small}\input{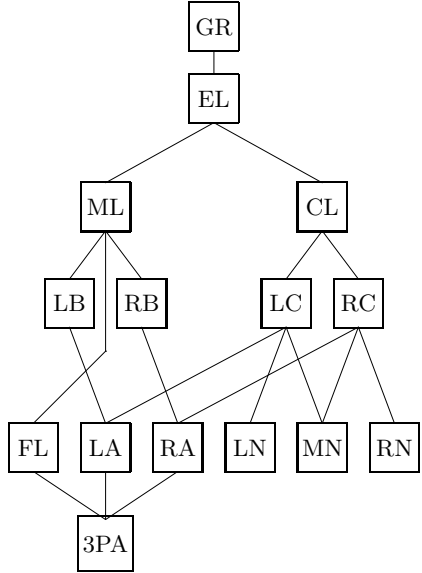}\end{small}
\caption{Varieties of loops of Bol-Moufang type and $3$-power associative
loops. If $A$, $B$ are varieties such that $A\subseteq B$ then $A$ is
depicted above $B$.}\label{Fg:Tree}
\end{figure}

\noindent We now show how the 14 varieties of loops of Bol-Moufang type are
related to each other.

\begin{lemma} The following inclusions hold among the varieties of loops of
Bol-Moufang type and \emph{3PA}: \emph{GR}$\subseteq$\emph{EL},
\emph{EL}$\subseteq$\emph{ML}, \emph{EL}$\subseteq$\emph{CL},
\emph{ML}$\subseteq$\emph{LB}, \emph{ML}$\subseteq$\emph{RB},
\emph{CL}$\subseteq$\emph{LC}, \emph{CL}$\subseteq$\emph{RC},\linebreak
\emph{ML}$\subseteq$\emph{FL}, \emph{LB}$\subseteq$\emph{LA},
\emph{RB}$\subseteq$\emph{RA}, \emph{LC}$\subseteq$\emph{LA},
\emph{RC}$\subseteq$\emph{RA}, \emph{LC}$\subseteq$\emph{LN},
\emph{LC}$\subseteq$\emph{MN}, \emph{RC}$\subseteq$\emph{MN},
\emph{RC}$\subseteq$\emph{RN}, \emph{FL}$\subseteq$\emph{3PA},
\emph{LA}$\subseteq$\emph{3PA}, \emph{RA}$\subseteq$\emph{3PA}. The situation
is depicted in Figure \emph{\ref{Fg:Tree}}.
\end{lemma}
\begin{proof} GR is contained in any variety of loops listed in Table \ref{Tb:Definitions}.
It is shown in \cite[Thm.\ 2]{Fe2} and \cite[Corollary 2]{CR} that extra
loops are precisely Moufang loops where every square belongs to the nucleus.
In \cite{Fe2}, Fenyves shows that extra loops are C-loops \cite{Fe2}, and
that C-loops are both LC-loops and RC-loops \cite[Thm.\ 4]{Fe2}. It is
well-known (cf \cite[Thm.\ 2.7]{Ro}) that Moufang loops are both left Bol and
right Bol. Moufang loops are flexible, as one can see upon letting $z=e$ in
$B15$. Robinson \cite[Thm.\ 2.1]{Ro} makes the simple observation that right
Bol loops are right alternative. The dual of this statement then holds, too.
LC-loops are left alternative by Lemma \ref{Lm:LC1}, and, dually, RC-loops
are right alternative. With the left alternative law at our disposal, we see
immediately that LC-loops are left nuclear square. The dual of this statement
then holds, too. Clearly, any of FL, LA or RA implies $3$-power
associativity. It remains to show that both LC-loops and RC-loops are middle
nuclear square. This follows from Lemma \ref{Lm:LC1}(iv) and its dual.
\end{proof}

\section{Distinguishing examples}\label{Sc:DE}

\noindent We proceed to show that all the 14 varieties of loops of Bol-Moufang
type are indeed distinct, and that no edges (inclusions) are missing in Figure
\ref{Fg:Tree}. Our intention is to come up with as few examples as possible to
accomplish this. It turns out that 8 examples and their duals suffice.

For the convenience of the reader, we provide Table \ref{Tb:Examples} that
points to examples distinguishing any two given varieties of loops of
Bol-Moufang type. If the cell in row A and column B of Table
\ref{Tb:Examples} is empty then A is a subvariety of B. If the cell contains
the integer $n$, then the loop of Example \ref{Sc:DE}.$n$ belongs to
$A\setminus B$. A primed number $n'$ indicates that one should use the dual
of the respective example.

\begin{table}\caption{Distinguishing varieties of loops of Bol-Moufang
type.}\label{Tb:Examples}
\begin{displaymath}
\begin{array}{c|cccccccccccccc}
    &\text{GR}&\text{EL}&\text{ML}&\text{CL}&\text{LB}&\text{RB}&\text{LC}
    &\text{RC}&\text{LA}&\text{FL}&\text{RA}&\text{LN}&\text{MN}&\text{RN}\\
    \hline
    \text{GR}&   &   &   &   &   &   &   &   &   &   &   &   &   &   \\
    \text{EL}&  1&   &   &   &   &   &   &   &   &   &   &   &   &   \\
    \text{ML}&  2&  2&   &  2&   &   &  2&  2&   &   &   &  2&  2& 2'\\
    \text{CL}&  3&  3&  3&   &  3& 3'&   &   &   &  3&   &   &   &   \\
    \text{LB}&  2&  2&  4&  2&   &  4&  2&  2&   &  4&  4&  2&  2& 2'\\
    \text{RB}&  2&  2& 4'&  2& 4'&   &  2&  2& 4'& 4'&   &  2&  2& 2'\\
    \text{LC}&  3&  3&  3&  5&  3& 3'&   &  5&   &  3&  5&   &   &  5\\
    \text{RC}&  3&  3&  3& 5'&  3& 3'& 5'&   & 5'&  3&   & 5'&   &   \\
    \text{LA}&  2&  2&  3&  2&  3& 3'&  2&  2&   &  4&  5&  2&  2& 2'\\
    \text{FL}&  2&  2&  6&  2&  6& 6'&  2&  2&  6&   & 6'&  2&  2& 2'\\
    \text{RA}&  2&  2&  3&  2&  3& 3'&  2&  2& 5'& 4'&   &  2&  2& 2'\\
    \text{LN}&  3&  3&  3&  7&  3& 3'&  7&  7&  7&  3&  7&   &  7&  5\\
    \text{MN}&  3&  3&  3&  8&  3& 3'&  8&  8&  8&  3&  8& 5'&   &  5\\
    \text{RN}&  3&  3&  3& 7'&  3& 3'& 7'& 7'& 7'&  3& 7'& 5'& 7'&
\end{array}
\end{displaymath}
\end{table}

All multiplication tables below have $0$ as a neutral element. We believe
that all examples below are as small as possible (when all properties are to
be satisfied at the same time).

\begin{example}[Extra loop that is not a group]
The Moufang loop that Goodaire et al. \cite{Goodaire} call 16/1 is a
nonassociative extra loop. Instead of giving its multiplication table, we
recall a general construction due to Chein \cite{Chein} that produces the loop
16/1.

For a group $G$, let $M(G,2)=G\times\{0,1\}$, where $(g,0)(h,0)=(gh,0)$,
$(g,0)(h,1)=(hg,1)$, $(g,1)(h,0)=(gh^{-1},1)$, and $(g,1)(h,1)=(h^{-1}g,0)$.
Then $M(G,2)$ is a nonassociative Moufang loop if and only if $G$ is
nonabelian.

Then $16/1$ is the loop $M(D_4,2)$, where $D_4$ is the dihedral group of order
$8$.
\end{example}

\begin{example}[Moufang loop that is neither left nuclear square nor middle
nuclear square] Take the loop $M(S_3,2)$, where $S_3$ is the symmetric group on
$3$ points.
\end{example}

\begin{example}[C-loop that is neither flexible, nor left Bol]
This example first appeared in \cite{KKP}.
\begin{displaymath}
\begin{array}{cccccccccccc}
 0& 1& 2& 3& 4& 5& 6& 7& 8& 9&10&11\\
 1& 2& 0& 4& 5& 3& 7& 8& 6&10&11& 9\\
 2& 0& 1& 5& 3& 4& 8& 6& 7&11& 9&10\\
 3& 4& 5& 0& 1& 2&10&11& 9& 8& 6& 7\\
 4& 5& 3& 1& 2& 0&11& 9&10& 6& 7& 8\\
 5& 3& 4& 2& 0& 1& 9&10&11& 7& 8& 6\\
 6& 7& 8&11& 9&10& 0& 1& 2& 4& 5& 3\\
 7& 8& 6& 9&10&11& 1& 2& 0& 5& 3& 4\\
 8& 6& 7&10&11& 9& 2& 0& 1& 3& 4& 5\\
 9&10&11& 7& 8& 6& 5& 3& 4& 0& 1& 2\\
10&11& 9& 8& 6& 7& 3& 4& 5& 1& 2& 0\\
11& 9&10& 6& 7& 8& 4& 5& 3& 2& 0& 1
\end{array}
\end{displaymath}
The loop is not flexible since $8(9\cdot 8)\ne (8\cdot 9)8$, and it is not left
Bol since $5(8(5\cdot 5))\ne (5(8\cdot 5))5$.
\end{example}

\begin{example}[Left Bol loop that is neither flexible, nor right alternative]
\begin{displaymath}
\begin{array}{cccccccc}
0& 1& 2& 3& 4& 5& 6& 7\\
1& 0& 3& 2& 5& 4& 7& 6\\
2& 3& 0& 1& 6& 7& 4& 5\\
3& 5& 1& 7& 0& 6& 2& 4\\
4& 2& 6& 0& 7& 1& 5& 3\\
5& 4& 7& 6& 1& 0& 3& 2\\
6& 7& 4& 5& 2& 3& 0& 1\\
7& 6& 5& 4& 3& 2& 1& 0
\end{array}
\end{displaymath}
The loop is not flexible since $1(2\cdot 1)\ne(1\cdot 2)1$, and it is not right
alternative since $6(4\cdot 4)\ne (6\cdot 4)4$.
\end{example}

\begin{example}[LC-loop that is neither right nuclear square, nor right alternative]
\begin{displaymath}
\begin{array}{cccccccccccc}
0& 1& 2& 3& 4& 5& 6& 7& 8& 9&10&11\\
1& 0& 3& 2& 5& 4& 7& 6& 9& 8&11&10\\
2& 3& 0& 1& 6& 7& 4& 5&10&11& 8& 9\\
3& 2& 6& 7& 0& 1&10&11& 4& 5& 9& 8\\
4& 5& 1& 0& 8& 9& 2& 3&11&10& 6& 7\\
5& 4& 8& 9& 1& 0&11&10& 2& 3& 7& 6\\
6& 7&11&10& 2& 3& 8& 9& 1& 0& 4& 5\\
7& 6&10&11& 3& 2& 9& 8& 0& 1& 5& 4\\
8& 9& 5& 4&11&10& 1& 0& 7& 6& 2& 3\\
9& 8& 4& 5&10&11& 0& 1& 6& 7& 3& 2\\
10&11& 7& 6& 9& 8& 3& 2& 5& 4& 0& 1\\
11&10& 9& 8& 7& 6& 5& 4& 3& 2& 1& 0
\end{array}
\end{displaymath}
The loop is not right nuclear square since $1(2(3\cdot 3))\ne (1\cdot 2)(3\cdot
3)$, and it is not right alternative since $1(2\cdot 2)\ne (1\cdot 2)2$.
\end{example}

\begin{example}[Flexible loop that is not left alternative]
\begin{displaymath}
\begin{array}{ccccc}
0& 1& 2& 3& 4\\
1& 0& 3& 4& 2\\
2& 4& 0& 1& 3\\
3& 2& 4& 0& 1\\
4& 3& 1& 2& 0
\end{array}
\end{displaymath}
The loop is not left alternative since $1(1\cdot 2)\ne (1\cdot 1)2$.
\end{example}

\begin{example}[Left nuclear square loop that is neither middle nuclear square,
nor $3$-power associative]
\begin{displaymath}
\begin{array}{cccccc}
0& 1& 2& 3& 4& 5\\
1& 5& 0& 4& 3& 2\\
2& 0& 4& 5& 1& 3\\
3& 4& 5& 0& 2& 1\\
4& 2& 3& 1& 5& 0\\
5& 3& 1& 2& 0& 4
\end{array}
\end{displaymath}
The loop is not middle nuclear square since $1((2\cdot 2)3)\ne (1(2\cdot 2))3$,
and it is not $3$-power associative since $1(1\cdot 1)\ne (1\cdot 1)1$.
\end{example}

\begin{example}[Middle nuclear square loop that is not $3$-power associative]
\begin{displaymath}
\begin{array}{cccccc}
0& 1& 2& 3& 4& 5\\
1& 2& 3& 0& 5& 4\\
2& 4& 5& 1& 3& 0\\
3& 5& 4& 2& 0& 1\\
4& 0& 1& 5& 2& 3\\
5& 3& 0& 4& 1& 2
\end{array}
\end{displaymath}
The loop is not $3$-power associative since $1(1\cdot 1)\ne(1\cdot 1)1$.
\end{example}

\section{Summary}

\noindent There are $14$ varieties of loops of Bol-Moufang type. Their
definitions can be found in Table \ref{Tb:Definitions}. They are related
according to Figure \ref{Fg:Tree}. One can look up examples distinguishing any
two varieties in Table \ref{Tb:Examples}. Since we believe this paper will be
used as a quick reference, we also include Table \ref{Tb:Ids}, that determines
the variety defined by any of the equations $Xij$ in $\mathcal B$, although the
same information is given in Section \ref{Sc:E}. To save space, we list
identities of type A, C, E as $Xij$ with $i<j$, and identities of type B, D, F
as $Xij$ with $i>j$.

\begin{table}\caption{Loop varieties determined by identities of Bol-Moufang type.}
\label{Tb:Ids}
\begin{small}
\setlength\arraycolsep{2pt}
\begin{displaymath}
\begin{array}{ccc}
\begin{array}{c|c|c|c|c|c}
    B\ldiv A&1&2&3&4&5\\
    \hline
    1&  &GR&LA&LC&LC\\
    2&GR&  &GR&GR&GR\\
    3&GR&EL&  &LC&LN\\
    4&LB&GR&GR&  &LA\\
    5&ML&GR&GR&FL&
\end{array}\quad\quad
&
\begin{array}{c|c|c|c|c|c}
    D\ldiv C&1&2&3&4&5\\
    \hline
    1&  &LA&GR&LC&CL\\
    2&GR&  &GR&MN&RC\\
    3&GR&ML&  &GR&GR\\
    4&GR&FL&ML&  &RA\\
    5&EL&GR&GR&GR&
\end{array}\quad\quad
&
\begin{array}{c|c|c|c|c|c}
    F\ldiv E&1&2&3&4&5\\
    \hline
    1&  &FL&GR&GR&ML\\
    2&RA&  &GR&GR&RB\\
    3&RN&RC&  &EL&GR\\
    4&GR&GR&GR&  &GR\\
    5&RC&RC&RA&GR&
\end{array}
\end{array}
\end{displaymath}
\end{small}

\end{table}

\section{Acknowledgement and historical remarks}

\noindent The classification of varieties of loops of Bol-Moufang type was
initiated by Fenyves \cite{Fe1}, \cite{Fe2}. He was aware of all results of
Section \ref{Sc:E} with the exception of the fact that C14 was equivalent to
the LC-identity A34, of the dual statement, and of some parts of Lemma
\ref{Lm:LC2}. He mentions all inclusions of Figure \ref{Fg:Tree} with the
exception of LC$\subseteq$MN and, dually, RC$\subseteq$MN. He only provides a
few distinguishing examples.

In the introduction of \cite{Fe2}, Fenyves claims: \emph{``Our results make
possible to decide of any two identities of Bol-Moufang type whether one of
them imply the other or not.''} This statement has now been justified.

The systematic notation is ours, and makes the discussion more transparent, in
our opinion.

Our investigations were aided by the equational reasoning tool \texttt{Otter}
and by the finite model builder \texttt{Mace4}. Both of these tools were
developed by McCune \cite{Mc}. Nevertheless, all proofs needed for the
classification (including those we only refer to) are now presented in full,
without any usage of computers.

We would like to thank the referee for several useful comments that allowed
us to substantially shorten the proofs of Lemmas \ref{Lm:LC1} and
\ref{Lm:LC2}.

%%%%%%%%%%%%%%%%%%%%%%%%%%%%%%%%%%%%%%%%%%%%%%%%%%%%%%%%%%%%%%%%%%%%%%%%%%%%%%%
% BIBLIOGRAPHY                                                                %
%%%%%%%%%%%%%%%%%%%%%%%%%%%%%%%%%%%%%%%%%%%%%%%%%%%%%%%%%%%%%%%%%%%%%%%%%%%%%%%
\bibliographystyle{plain}

\end{document}